\documentstyle[12pt]{article}
\newcommand{\const}{\mathop{\rm const}\limits}

\newcommand{\Law}{\mathop{\rm Law}\limits}

\newcommand{\Var}{\mathop{\rm Var}\limits}

\newcommand{\Ent}{\mathop{\rm Ent}\limits}

\newcommand{\Cov}{\mathop{\rm Cov}\limits}

\begin{document}

\begin{center}

{\bf UNBIASED MONTE CARLO ESTIMATION  FOR SOLVING  }\\

\vspace{4mm}

{\bf OF LINEAR INTEGRAL EQUATION,  }\par

\vspace{4mm}

{\bf with error estimate. }\par

\vspace{4mm}

 $ {\bf E.Ostrovsky^a, \ \ L.Sirota^b } $ \\

\vspace{4mm}

$ ^a $ Corresponding Author. Department of Mathematics and computer science, Bar-Ilan University, 84105, Ramat Gan, Israel.\\
\end{center}
E-mail: \ eugeneiostrovsky@gmail.com \  eugostrovsky@list.ru\\

\begin{center}

$ ^b $  Department of Mathematics and computer science. Bar-Ilan University,
84105, Ramat Gan, Israel.\\

E-mail: \ sirota3@bezeqint.net \\

\vspace{4mm}
                    {\sc Abstract.}\\

 \end{center}

 \vspace{4mm}

   We  offer a new  Monte-Carlo method for solving linear integral equation which gives the unbiased estimation for solution
of Volterra's and Fredholm's type, and  consider the problem of confidence region building.\par
  We study especially  the case of the so-called  equations with weak singularity in the kernel of Abelian type. \par

\vspace{4mm}


\vspace{4mm}

 \par

\vspace{4mm}

\section{Notations. Statement of problem.}

\vspace{3mm}

  We intent in this article to study the numerical Monte-Carlo method for solving of  the linear integral
 equation, for example, of a form

$$
x(t) = f(t) + \lambda \int_0^t K(t,s) \ x(s) \ ds, \eqno(1.0)
$$
Volterra's equation. (The case of Fredholm's equation will be considered further.) \par
 Here $  x(t) $ is unknown function, $ K(t,s) $ is kernel, the function $  f(t) $ is "right-hand"  side,
$ \lambda $ is {\it  positive } number.\par
 Briefly:

$$
x = f + \lambda \ K [x], \eqno(1.0a)
$$
where $ K[x] = K[x](t) $ is linear integral operator of Volterra's type:

$$
K[x](t) = \int_0^t K(t,s) \ x(s) \ ds.
$$

\vspace{3mm}

 {\bf The offered here method  gives as ordinary the optimal rate of convergence, exponential tail estimation  for confidence probability,
but in addition our (random!) estimates of solution are unbiased. } \par

\vspace{3mm}

  We will consider the equation (1.0), as well as all the next equations, in the space of continuous functions $ C(T)  $ defined on the set
 $  t \in [0,T],  $  where $ T = \const \in (0, \infty), $  suppose therefore $ f(\cdot) \in C(T), \ K(\cdot, \cdot) \in C(T \times T), $
 and denote as usually

$$
||f|| = \max_{t \in [0,T]} |f(t)|, \hspace{6mm} ||K|| = \max_{s,t \in [0,T]} ||K(t,s)||. \eqno(1.1)
$$

\vspace{3mm}

 These equations appear in particular in the reliability theory: renewal equation, equating of periodical checking etc., see, e.g.
the article \cite{Grigorjeva1}, where  is explained in particular why the Monte-Carlo method is natural for this problem. \par

\vspace{3mm}

 It is well known that this problem is well posed: the solution $  x(t) $ there exists, is unique and dependent continuously on the
entries:  $ f(\cdot), \ K(\cdot, \cdot)  $  and may be computed by means of the standard recursion procedure:

$$
x_n = f + \lambda \  K \left[x_{n-1} \right],  \hspace{6mm}, x_0 = 0, \  x_1 = f; n = 2,3,\ldots \eqno(1.2)
$$

\vspace{4mm}

\section{ Essence of offered method.  }

\vspace{3mm}

  The solution $ x(\cdot) $ may be represented by means of the so - called  Neuman's series:

$$
x = f + \sum_{n=1}^{\infty} \ \lambda^n K^n[f] = \sum_{n=0}^{\infty} \ \lambda^n K^n[f], \eqno(2.1)
$$
where $ K^n $ denotes the $ n^{th} $ iteration (power) of the operator $  K, $ so that $ K^0[f] = f; $ which converges uniformly.\par

 Denote also $ y_n(t) = K^n[f](t);  $ then $ x = \sum_{n=0}^{\infty} y_n $ and we have for the values $ n = 1,2,\ldots: \
y_n(t) = \lambda^n \times $

$$
\int_0^t ds_1 \int_0^{s_1} ds_2 \ldots \int_0^{s_{n-1}} ds_n \cdot K(t,s_1) K(s_1,s_2) \ldots K(s_{n-1}, s_n) \ f(s_n) =
 \lambda^n  t^n \times
$$

$$
\int_0^1 ds_1 \int_0^{s_2} \ldots \int_0^{s_{n-1}} ds_n \cdot K(t,t s_1) K(ts_1,ts_2) \ldots K(ts_{n-1}, ts_n) \ f(ts_n).
 \eqno(2.2)
$$

 Let us introduce the following $ n \ - $ dimensional simplex (polygon) $  S(n) :=  $

$$
 \{s_1,s_2,\ldots,s_n: 0 < s_1 < 1, \ 0 < s_2 < s_1, 0 < s_3 < s_2, \ldots, 0 < s_n < s_{n-1} \}, \eqno(2.3)
$$
and denote $ L_n(t, \vec{s}) = L_{n;K,f}(t, \vec{s}) =  L_n(t, s)= $

$$
 K(t,t s_1) K(ts_1,ts_2) \ldots K(ts_{n-1}, ts_n) \ f(ts_n), \ \vec{s} = s = (s_1, s_2, \ldots,s_n);
$$
then

$$
y_n(t) = \int_{S(n)} L_n(t, \vec{s}) \ ds. \eqno(2.4)
$$
 Note that the volume of the polygon $  S(n) $ is equal to $ 1/n!, $ and if we introduce the {\it  probability } measure $ \mu_n(\cdot) $
on the Borelian subsets of the simplex $  S(n) $ with a density $ n! \ ds: $

$$
\mu_n(A) = n! \int_A  ds,
$$
then the expression for $  x(t) $ may be rewritten as follows:

$$
x(t) = f(t) + \sum_{n=1}^{\infty} \frac{\lambda^n \ t^n}{n!} \int_{S(n)} L_n(t, \vec{s}) \ \mu_n(ds).\eqno(2.5)
$$

 If we denote $  x_{\lambda} = x_{\lambda}(t) =  e^{- \lambda t} \cdot x(t),  $ then

$$
x_{\lambda}(t) = e^{-\lambda t} f(t) + \sum_{n=1}^{\infty}  e^{-\lambda t} \frac{\lambda^n \ t^n}{n!} \int_{S(n)} L_n(t, \vec{s}) \ \mu_n(ds), \eqno(2.6)
$$
which may be formally rewritten as follows:

$$
x_{\lambda}(t) = \sum_{n=0}^{\infty}  e^{-\lambda t} \frac{\lambda^n \ t^n}{n!} \int_{S(n)} L_n(t, \vec{s}) \ \mu_n(ds), \eqno(2.6a)
$$

 Let us introduce a sufficiently rich probability space $  (\Omega, B,{\bf P}) $ with probability $ {\bf P,}  $
 expectation  $ {\bf E}  $ and variance $ \Var; $ and the r.v. $  \tau $ which has a {\it Poisson  } distribution with parameter $ \lambda \cdot t $
(which is closely related with Poisson flow of intensity $ \lambda): $

$$
{\bf P} (\tau = n) = e^{-\lambda t} \frac{\lambda^n \ t^n}{n!}; \eqno(2.7)
$$
then the expression  (2.6a) takes the form

$$
x_{\lambda}(t) = {\bf E} \int_{S(\tau)} L_{\tau}(t, \vec{s}) \ \mu_{\tau}(ds), \eqno(2.6a)
$$

 Further, the integral in the right - hand side (2.6a) may be represented as follows. We introduce the random
vector $ \xi_{\tau} $  of a dimension $ \tau  $ which has the uniform (conditional) distribution in the simplex $  S(\tau): $

$$
{\bf P} ( \xi_{\tau} \in A  )/\tau = \mu_{\tau}(A), \eqno(2.7)
$$
then

$$
\int_{S(\tau)} L_{\tau}(t, \vec{s}) \ \mu_{\tau}(ds) = {\bf E}_{\tau} L_{\tau}(t, \vec{\xi_{\tau}}), \eqno(2.8)
$$

so that

$$
x_{\lambda}(t) = {\bf E}  L_{\tau}(t, \vec{\xi_{\tau}}). \eqno(2.9)
$$

 By the {\it fixed} value $  \tau $ the integral in the expression (2.8) may be computed by means of the Monte Carlo method:

$$
\int_{S(\tau)} L_{\tau}(t, \vec{s}) \ \mu_{\tau}(ds) \approx \frac{1}{k} \sum_{j=1}^k  L_{\tau}(t, \vec{\xi}^{(j)}_{\tau}), \eqno(2.10)
$$
where $ \vec{\xi}^{(j)}_{\tau} $ are independent copies of the r.v. $ \vec{\xi}_{\tau}. $ \par

 \vspace{4mm}

{\it Correspondingly, the Monte Carlo  approximation  for the whole value $ x_{\lambda}(t) $ may be offered as follows:}

$$
\hat{x}_{\lambda} =\hat{x}_{\lambda, N,Z}(t)  \stackrel{def}{=} \frac{1}{Z} \sum_{i=1}^Z \
 \left[ \frac{1}{N(\tau(i))} \sum_{j=1}^{ N(\tau(i)) }  L_{\tau(i)} \left(t, \vec{\xi}^{(j)}_{\tau(i)} \right) \right], \eqno(2.11)
$$
{\it where  $  \tau(i) $ are independent copies of the r.v. $ \tau, $ i.e. are independent Poisson distributed r.v. with parameter $ \Lambda = \lambda \ t,  $
and $  N = N(n) $ be some non-random positive numerical sequence, her choice will be specified later. } \par

\vspace{4mm}

 We state by definition

$$
\frac{1}{N(n)}\sum_{j=1}^{N(n)} L_j \stackrel{def}{=} 0, \eqno(2.11a)
$$
if $  N(n) = 0. $ \par

\vspace{3mm}

 Evidently, the approximation $  \hat{x}_{\lambda} =\hat{x}_{\lambda, N,Z}(t) $ of the solution $ x_{\lambda}(t) $ is unbiased:
$ {\bf E}\hat{x}_{\lambda, N,Z}(t) = x_{\lambda}(t). $  \par

 Another approach for Monte Carlo solving of linear integral equation which gives  biased estimation see in the article
 \cite{Ostrovsky302}. \par

\vspace{3mm}

Let's count the amount $  R   $ of elapsed r.v. for the $ \hat{x}_{\lambda} =\hat{x}_{\lambda, N,Z}(t) $ computation.

$$
R = Z \cdot \sum_{i=1}^Z \tau(i) N(\tau(i)), \eqno(2.12)
$$
so that $ R $ is {\it  random variable } with the expectation

$$
\Theta \stackrel{def}{=} {\bf E} R \asymp Z \cdot \sum_{n=1}^{\infty} e^{- \Lambda} \frac{\Lambda^n}{n!} \cdot (n N(n)), \eqno(2.13)
$$
where $  \Lambda = \lambda t. $\par

\section{ Error estimate. }

\vspace{3mm}

 Let us estimate the variance of the approximate solution $ \hat{x}_{\lambda} =\hat{x}_{\lambda, N,Z}(t). $ Note first of all that

$$
\Var \left[\hat{x}_{\lambda} \right] = \frac{1}{Z} \cdot
\Var \left[ \frac{1}{N(\tau)} \sum_{j=1}^{N(\tau)} L_{\tau} \left(t, \vec{\xi}^{(j)}_{\tau} \right) \right].
\eqno(3.1)
$$

 Further, we will use the next formula

$$
\Var(\zeta) =  {\bf E} \left\{ {\bf  E  } ( \zeta - {\bf E}\zeta )^2 /\tau \right\}:
$$

$$
Z \cdot \Var \left[\hat{x}_{\lambda} \right] \le ||f||^2 \cdot e^{-\Lambda} \sum_n \frac{\Lambda^n \ ||K||^{2n}}{N(n) \ n!^2}=
$$

$$
||f||^2 \cdot e^{-\Lambda} \sum_n \frac{Q^n}{N(n) \ n!^2}, \hspace{7mm}  Q := \Lambda \ ||K||^2.  \eqno(3.2)
$$

\vspace{3mm}

{\bf  Lemma 3.1. }  Let us consider the following constrained optimization problem:

$$
D := \min_{  \{N(n)\} \ge 1 } \sum_n \frac{A(n)}{N(n)} \ / \ \left\{ \sum_n B(n) N(n) \le M \right\}. \eqno(3.3)
$$
 Here $ \{ A(n)  \}, \  \{ B(n) \} $ are (may be unbounded) positive sequences, $ M = \const >> 1,
 \ N(n) >  0. $ \par

 We derive using the Lagrange's factors method that

$$
D = \frac{ \left[ \sum_n \sqrt{ A(n) \ B(n)} \right]^2}{M}, \eqno(3.4)
$$
and is achieved iff

$$
N(n) := N_0(n) \stackrel{def}{=} \frac{M}{\sum_n \sqrt{ A(n) \ B(n) }}  \cdot \sqrt{ \frac{A(n)}{B(n)}}, \eqno(3.5)
$$
up to rounding to the nearest greatest integer number:

$$
N(n) := \Ent[ N_0(n)  ] + 1,
$$
where  $ \Ent(z) $  denotes the integer part of real positive number $ z, $ if it is known that the numbers  $  N(n) $
are integer. \par

 In the last case we have only

$$
D \asymp \frac{ \left\{ \sum_n \sqrt{ A(n) \ B(n)} \right\}^2}{M}. \eqno(3.6)
$$

\vspace{3mm}

 Let $ \Theta = {\bf E}R $ be a fixed "great" number, for instance, $  10^7 - 10^8. $ It seems the following {\it constrained}
optimization problem of the variance minimization:

$$
\sum_n \frac{Q^n}{N(n) \ n!^2} \to \min / \sum_{n=1}^{\infty} \frac{\Lambda^n}{n!} \cdot (n N(n)) = e^{\Lambda} \Theta/Z.
$$
 Of course, it will be supposed that $ \Theta \ e^{\Lambda} >> Z. $ \par

We observe  using lemma 3.1

$$
\min_{ \{ N(n) \} } \left[ \frac{ \Var( Z \cdot \hat{ x_{\lambda}} )}{ e^{-\Lambda}||f||^2}\right] \approx
\frac{1}{\Theta} \cdot \left( \sum_n \frac{\Lambda^n \ ||K||^n  \ \sqrt{n}}{(n!)^{3/2}} \right)^2, \eqno(3.7)
$$
wherein the (quasi \ - ) optimal values $ N(n) = N_0(n) $ are following:

$$
N_0(n)= \Ent \left[ \frac{e^{\Lambda} Q/Z}{ \left[ \sum_m  \Lambda^m \ ||K||^m \ (m!)^{-3/2} \ \sqrt{m} \right]  } \cdot ||K||^n \cdot
\frac{\sqrt{n}}{ {n!}} \right] + 1, \eqno(3.8)
$$
wherein by the practical using in the case when the value $ N_0(n) $ is sufficiently small, for instance, if $ N_0(n) \le 10, $
then we oblige to take $ N_0(n):=0; $  see (2.11a). \par

 This imply that the rate of convergence of offered method is equal to $  1/\sqrt{Z}, $ as in the classical Monte Carlo method.\par

 It remains  to use the classical Central Limit Theorem (CLT) to construct the confidence interval for the solution $  x_{\lambda}(t). $\par

\vspace{4mm}

\section{ Fredholm's equations.}

\vspace{3mm}

 We consider in this section the Monte - Carlo approach for solution of Fredholm's \cite{Fredholm1}  linear integral equation of a second kind

$$
x(u) = f(u) + \lambda \ \int_V K(u,v) \ x(v) \ \nu(dv), \eqno(4.1)
$$
or equally

$$
x(u) = f(u) + \lambda \  K[x](u), \eqno(4.1a)
$$
where $ \lambda = \const \in (0,1), \ (V, A, \nu) $ with a distance $ \rho = \rho(v_1, v_2) $ is compact metric measurable
probability space: $ \nu(V) = 1;   \ u,v \in V,  $ and both the functions $  f(\cdot), K(\cdot, \cdot) $ are continuous, and we
denote  as in the first section

$$
||f||= \max_{u \in V}|f(u)|, \ ||K|| = \max_{u,v \in V} |K(u,v)|.
$$

 Another (deterministic) approach via the so-called Fredholm's determinants \cite{Fredholm1} computing implementation
see in the article \cite{Bornemann1}. \par

 The norm of linear operator $  K[x] (u) = \int_V K(u,v) x(v)  \nu (dv) $ in the space of continuous functions $ C(V) = C(V, \rho) $
 will be denoted by $ |||K|||; $ it is known that

$$
|||K||| = \max_{u \in V} \int_V |K(u,v)| \ \nu(dv). \eqno(4.2)
$$
 Evidently, $ |||K||| \le ||K||. $ \par
Further, let us denote

$$
r_n = r_n(K) = ||| K^n|||^{1/n}; \hspace{6mm} r = r(K) = \lim_{n \to \infty} r_n(K). \eqno(4.3)
$$
 The last limit there exists and is named as spectral radius of the operator $  K; $ see \cite{Dunford1}, chapters 4,5. \par

\vspace{3mm}

{\it  We suppose at first that } $ |||K||| \le 1; $ then the continuous solution  $ x = x(u) $ of (4.1) there exists, is unique
and may be represented by means of uniform converge Neuman series:

$$
x(u) = f(u) + \sum_{n=1}^{\infty} \lambda^n \ K^n[f](u) = f(u) + \sum_{n=1}^{\infty} \lambda^n \  y_n(u), \hspace{5mm} y_n(u) := \eqno(4.4)
$$

$$
  \int_V \nu(ds_1) \int_V \nu(ds_2) \ldots \int_V \nu(ds_n) K(u,s_1)K(s_1,s_2) \ldots K(s_{n-1}, s_n) \ f(s_n). \eqno(4.5)
$$
 Recall that $  0 < \lambda < 1. $\par

 Let $ \gamma(i), \ i = 1,2,\ldots,n $ be independent random variables with distribution $ \nu:  $

$$
{\bf P}(\gamma(i) \in A) = \nu(A).
$$
 Then the function $ y_n(\cdot) $  has a probabilistic representation: $ y_n(u) = $

$$
 y_n(u) = {\bf E} K(u,\gamma(1)) \ K(\gamma(1), \gamma(2)) \ \ldots \ K(\gamma(n-1), \gamma(n)). \eqno(4.6)
$$

 If we denote (in this section)

$$
x_{\lambda}(u) = (1-\lambda) x(u), \ \vec{s}_n = \{ s(1), s(2), \ldots, s(n) \},
$$

$$
 L_n(u, \vec{s}_n) = K(u,s(1)) K(s(1), s(2)) \ldots K(s(n-1), s(n)),
$$

$$
\vec{\theta}_n = \theta_n = \{ \gamma(1), \gamma(2), \ldots, \gamma(n)  \},
$$
then

$$
y_n(u) = {\bf E} L_n(u, \vec{\theta}_n)
$$
and

$$
x_{\lambda} = \sum_{n=0}^{\infty} (1-\lambda) \ \lambda^n {\bf E}  L_n(u, \vec{\theta}_n). \eqno(4.7)
$$

 The non - negative integer valued random variable $ \Delta $ with distribution

$$
{\bf P} ( \Delta = n ) = (1 - \lambda) \ \lambda^n, \ n = 0,1,2, \ldots
$$
is named (integer) geometrical distributed, write: $ \Law(\Delta) = G_{\lambda}.  $ \par
 We can rewrite the expression (4.7) using this notation as follows

$$
x_{\lambda} = {\bf E}  L_{ \Delta }(u, \vec{\theta}_{\Delta}). \eqno(4.8)
$$

 The Monte - Carlo approximation  $  \hat{x}_{\lambda} = \hat{x}_{\lambda}(u) = \hat{x}_{\lambda; Z, \{ N(n) \}}(u)  $ for $ x_{\lambda} $
may be written as before

$$
\hat{x}_{\lambda}(u) = \frac{1}{Z} \sum_{i=1}^Z \frac{1}{N(\Delta(i))} \sum_{j=1}^{ N(\Delta(i)) } L_{\Delta(i)}(u, \theta^{(j)}_{\Delta(i)}), \eqno(4.9)
$$
where $ \{ \Delta(i)  \} $ are independent copies of the r.v. $  \Delta $ and  $ \theta^{(j)}_{\Delta(i)} $ are independent copies of the
random vector $ \theta_{\Delta(i)}. $\par

 In order to calculate (estimate) the variance of $ \hat{x}_{\lambda}(u), $ we need to use the following definition. Let $ K[x](u) $ be any linear
integral operator with kernel $  K: $

$$
K[x](u) = \int_V K(u,v) \ x(v) \ \nu(dv).
$$
 The {\it  Kroneker's  } square  $  K^{[2]}  $  of the operator $  K $  is an operator acting as follows:

$$
K^{[2]}[x](u) \stackrel{def}{=} \int_V K^2(u,v) \ x(v) \ \nu(dv). \eqno(4.10)
$$

{\it  We impose another condition on the coefficient $  \lambda $ and on the kernel} $ K: $

$$
\lambda \ \cdot |||K^{[2]}||| < 1. \eqno(4.11)
$$

 Then

$$
\frac{Z}{(1 - \lambda) ||f||^2} \Var \hat{x}_{\lambda} \le \sum_n \frac{\lambda^n |||K^{[2]}|||^n}{N(n)}. \eqno(4.12)
$$

\vspace{3mm}

 Let's count now the amount $  R_F   $ of elapsed {\it standard, i.e. uniform [0,1] distributed, }
 r.v. for the $ \hat{x}_{\lambda} =\hat{x}_{\lambda, N,Z}(t) $ computation.\par

 \vspace{4mm}

{\it  We suppose that for one $ \theta^{j}_k  $
random vector number  generation are  need exactly  $ k \cdot d, \ d = \const = 1,2,\ldots  $ standard r.v.} \par

\vspace{4mm}

 Then

$$
R_F = Z  \cdot \sum_{i=1}^Z \tau(i) N(\tau(i)), \eqno(4.13)
$$
so that $ R_F $ is {\it  random variable } with the expectation

$$
\Theta_F \stackrel{def}{=} {\bf E} R_F \asymp Z \cdot d \cdot \sum_{n=0}^{\infty} \left[ (1- \lambda) \ \lambda^n \cdot (n N(n)) \right]. \eqno(4.14)
$$

  We get denoting

$$
M_F = \frac{\Theta_F}{(1 - \lambda) \ Z \ d }
$$
to the following constrained optimization problem assuming $  M_F $ a fixed number of large:

$$
\left[\frac{(1 - \lambda)||f||^2}{Z} \right] \times
\sum_n \frac{\lambda^n |||K^{[2]}|||^n}{N(n)} \to \min / \left[ \sum_n \lambda^n \cdot (n N(n)) = M_F  \right]. \eqno(4.14)
$$

 Let us introduce the following function:

$$
G_{\alpha}(z) = \sum_{n=0}^{\infty} n^{\alpha} \ z^n; \ \alpha = \const \ge 0, \ 0 \le x < 1.
$$
 We find  tacking into account the proposition of lemma 3.1:

$$
D_F := \min_{ \{  N(n) \} } \Var \hat{x}_{\lambda}  \asymp \frac{(1 - \lambda)^2 \ d \ ||f||^2 }{\Theta_F} \cdot
 G_{1/2} \left( \lambda \sqrt{ |||K^{[2]}||| } \right) \eqno(4.15)
$$
and this minimum is attained iff

$$
N(n) := 1 + \Ent(N_0(n)),
$$
where

$$
N_0(n) \stackrel{def}{=}  \left[ \frac{M_F}{G_{1/2} \left( \lambda \sqrt{ |||K^{[2]}||| } \right)} \right] \times
\frac{|||K^{[2]}|||^{n/2}}{\sqrt{n}}. \eqno(4.16)
$$
 Of course,  if $ N_0(n) \le 10, $ we must take $  N(n) = 0. $ \par

 We conclude again that the rate of convergence offered estimate is equal to $ 1/\sqrt{ \Theta_F }, $ as in the one dimensional case.\par

\vspace{3mm}

{\bf Remark 4.1. } Note that as $  z \to 1-0 $

$$
G_{1/2}(z) \sim 0.5 \ \sqrt{\pi} \ |\ln z|^{-3/2}.
$$

\vspace{4mm}

\section{Equations with weak singularity.}

\vspace{3mm}

 We consider in this section the Volterra's type integral equation with weak (Abel's) singularity of a form

$$
x(t) = f(t) + \lambda \int_0^t \frac{K(t,s) \ x(s) \ ds}{(t-s)^{\alpha}}, \eqno(5.1)
$$
where as before $ \lambda = \const > 0, \ t \in [0.T], \ T = \const > 0, \  f(\cdot), \ K(\cdot, \cdot) $ are continuous functions,
$ \alpha := 1 - \beta = \const \in (0,1); $ the case $ \alpha = 0 $ was considered in the second section.\par

This case $ \alpha > 0 $ has a number of interesting features, and we will briefly enumerate.\par

 First of all  $ x(t) = f(t) + \sum_{n=1}^{\infty} y_n(t), $  where  $ y_n(t) = \lambda^n \times $

$$
\int_0^t ds_1 \int_0^{s_1} ds_2 \ldots \int_0^{s_{n-1}} ds_n \
 \frac{ K(t,s_1) K(s_1,s_2) \ldots K(s_{n-1}, s_n) \ f(s_n)}{ \left[(t-s_1)(s_1 - s_2) \ldots (s_{n-1} - s_n)  \right]^{\alpha} } =
 [\lambda t^{\beta}]^n \times
$$

$$
\int_0^1 ds_1 \int_0^{s_2} \ldots \int_0^{s_{n-1}} ds_n \
\frac{ K(t,t s_1) K(ts_1,ts_2) \ldots K(ts_{n-1}, ts_n) \ f(ts_n)}{ \left[ (1-s_1)(s_1 - s_2) \ldots (s_{n-1} - s_n)\right]^{\alpha}}.
 \eqno(5.2)
$$

 Denote

$$
R_{\alpha,n} ( \vec{s}  ) = (1-s_1)^{-\alpha} (s_2 - s_1)^{-\alpha} (s_3 - s_2)^{-\alpha} \ldots (s_n - s_{n-1})^{-\alpha},
$$

$$
W_n(\beta) = \frac{\Gamma^n(\beta)}{\Gamma(1 + n \beta)},
$$
where $ \Gamma(\cdot) $ is ordinary Gamma function. Evidently, $ \lim_{n \to \infty} W_n(\beta) = 0.  $\par
  Note that

$$
\int \int \ldots \int_{S(n)}\frac{ds_1 ds_2 \ldots ds_n}{(1-s_1)^{\alpha} (s_2 - s_1)^{\alpha} (s_3 - s_2)^{\alpha} \ldots (s_n - s_{n-1})^{\alpha} } =
W_n(\beta).
$$

\vspace{3mm}

{\it  The following function $  h_{\alpha}(s) =  h_{\vec{\alpha}}(s), \ s \in S(n),  $
  could be chosen as a density of distribution with support on the   simplex} $  S(n): $

$$
h_{\alpha}(s) = \frac{R_{\alpha,n}(s)}{W_{n}(\beta)}.
$$

\vspace{3mm}

{\bf Definition 5.1.} (See \cite{Ostrovsky302}.)  The random vector $  \kappa = \kappa_{\alpha,n} = \vec{\kappa}  = \vec{\kappa}_{\alpha,n}$
with values in the polygon $  S(n) $ has by definition a {\it polygonal Beta distribution,} write:   $  \Law(\kappa) = PB(\alpha,n), $   iff it has a
density $  h_{\alpha}(s), \ s \in S(n).  $\par

 On the other word,

$$
{\bf P}(\kappa \in G) = \int_G h_{\alpha}(s) \ ds \stackrel{def}{=} \mu_{\alpha,n}(G), \ G \subset S(n).
$$

 Evidently, $ \mu_{\alpha,n}(\cdot) $ is the probabilistic Borelian measure on the set $  S(n).  $ \par

 The expression for the function $ y_n(\cdot) $ may be rewritten as follows:

 $$
  y_n(t) = \lambda^n t^{n\beta} \cdot W_n(\beta) \cdot {\bf E} L_n(t \ \vec{\eta}_n ),\eqno(5.3)
 $$
where the random vector  $ \vec{\eta}_n $ has the polygonal  beta distribution of dimension $ n $ with parameters $ (\alpha,n).  $ \par

 Recall that the function of Mittag - Lefler $ E_{\beta}(z), $ more exactly, the family of the functions which dependent on the
positive real parameter $ \beta, \ \beta > 0, $  is defined for all (may be complex) values $  z  $ by the formula

$$
E_{\beta}(z) = \sum_{n=0}^{\infty} \frac{z^n}{\Gamma(1 + n \beta)}, \ \beta = \const > 0.
$$

 We define also some slight generalization of the Mittag - Lefler function:

$$
E_{\beta, \alpha, \delta}(z) = \sum_{n=1}^{\infty} \frac{z^n}{n^{\delta} \ \Gamma^{\alpha}(1 + n \beta)}, \ \alpha, \beta = \const > 0, \ \delta = \const.
$$
 Evidently, $ E_{\beta,1, 0}(z) = E_{\beta}(z)-1.  $ \par

 This definition with investigation of properties of this function belongs to G.Mittag-Lefler \cite{Mittag1}. See also a recent
article  \cite{Grothaus1}, (with reference therein,) where are described in particular some interest applications of these
functions.\par

\vspace{3mm}

{\bf Definition 5.2.} The {\it integer valued} non - negative random variable $  \zeta $ has by definition Mittag-Lefler distribution
with parameters $ (\beta, \mu), \ \beta, \mu > 0, $ write:  $ \Law(\zeta) = R_{\beta}(\mu), $ iff

$$
{\bf P} (\zeta = n) = \frac{\mu^n/\Gamma(1 + \beta n)}{E_{\beta}(\mu)}, \ n = 0,1,2,\ldots \eqno(5.4)
$$

\vspace{3mm}

{\bf Remark 5.1.} Our definition (5.2) is unlike from ones in the article \cite{Grothaus1}, where was introduced the so - called
{\it   continuous } Mittag - Lefler distribution. \par

\vspace{3mm}

 Denote

 $$
 x_{\lambda, \beta}(t) = \frac{x(t)}{E_{\beta}(\Lambda_{\beta}(t))}.
 $$

 The function $  x_{\lambda, \beta}(t), $ which is proportional to the true
  solution $ x(t) $ of the equation (5.1), may be probability represented  as follows

$$
x_{\lambda, \beta}(t) = {\bf E}  L_{\zeta}(t, \vec{\eta_{\zeta}}). \eqno(5.5)
$$

 By the {\it fixed} value $  \zeta $ the integral in the expression (5.5) may be computed by means of the Monte Carlo method:

$$
\int_{S(\zeta)} L_{\zeta}(t, \vec{s}) \ \mu_{\zeta}(ds) \approx \frac{1}{k} \sum_{j=1}^k  L_{\zeta}(t, \vec{\eta}^{(j)}_{\zeta}),
$$
where $ \vec{\eta}^{(j)}_{\zeta} $ are independent copies of the r.v. $ \vec{\eta}_{\zeta}. $ \par

 Correspondingly, the Monte Carlo  approximation  for the whole value $ x(t) $ may be offered as in the second section as follows:

$$
\hat{x}_{\lambda, \beta} =\hat{x}_{\lambda, \beta; N,Z}(t)  \stackrel{def}{=} \frac{1}{Z} \sum_{i=1}^Z \
 \left[ \frac{1}{N(\zeta(i))} \sum_{j=1}^{ N(\zeta(i)) }  L_{\zeta(i)} \left(t, \vec{\eta}^{(j)}_{\zeta(i)} \right) \right], \eqno(5.6)
$$
where  $  \zeta(i) $ are independent copies of the r.v. $ \zeta, $ i.e. are independent Mittag - Lefler's distributed r.v. with parameter
$ \beta, \Lambda_{\beta} = \Lambda_{\beta}(t) \stackrel{def}{=} \lambda \ t^{\beta},  $ and $  N = N(n) $ be some non-random {\it integer}
positive numerical sequence, her choice will be specified later. \par

 Further, let us estimate the variance of our estimation $ \hat{x}: $

$$
Z \ \Var{\hat{x}_{\lambda, \beta}} \le \sum_n W_n^2(\beta) \frac{\Lambda^{n}_{\beta}(t)}{N(n)}\int_{S(n)}L^2(t,s) \ \mu_n(ds) \le
$$

$$
||f||^2 \ \sum_n  \frac{W^2_n(\beta) \ \Lambda^{n}_{\beta}(t) \ ||K||^{2n} }{ N(n)  }.\eqno(5.7)
$$
 Let us discuss now the question of amount $  R_{\beta}  $ elapsed random variables.  We find analogously to the second section

$$
R_{\beta} = Z \cdot \sum_{i=1}^Z \zeta(i) N(\zeta(i)), \eqno(5.8)
$$
so that $ R_{\beta} $ is {\it  random variable } with the expectation

$$
\Theta_{\beta} \stackrel{def}{=} {\bf E} R_{\beta} \asymp Z \cdot \
\sum_{n=1}^{\infty}\frac{ \Gamma^n(\beta) \ \Lambda_{\beta}^n/\Gamma(1 + \beta n)}{E_{\beta}(\Lambda_{\beta})} \cdot ( n \ N(n)) . \eqno(5.9)
$$

\vspace{4mm}

 Let $ \Theta_{\beta} = {\bf E}R $ be a fixed "great" number. Denote

$$
M_{\beta} =  E_{\beta}( \Lambda_{\beta}(t) ) \cdot \frac{\Theta_{\beta}}{Z};
$$
it will be presumed that $  M_{\beta} $  is also a great number,  for instance, $  10^7 - 10^8. $ \par

\vspace{3mm}

 It seems the following {\it constrained} optimization problem of the variance minimization:

$$
||f||^2 \sum_n \frac{ W_n^2(\beta) \ ||K||^{2n} \ \Lambda^n_{\beta}}{N(n)} \to \min /
\sum_{n} \frac{\Lambda^n_{\beta} \ \Gamma^n(\beta)}{\Gamma(1 + \beta n)} \cdot (n N(n)) \le M_{\beta}. \eqno(5.10)
$$

We deduce  using again lemma 3.1

$$
\min_{ \{ N(n) \} }   \Var \hat{ x_{\lambda,\beta}} \approx \frac{||f||^2}{ M_{\beta} } \cdot
 E^2_{ \beta, \ 3/2, \ 1/2} \left( \Lambda_{\beta} \ ||K|| \  \Gamma^{3/2}(\beta)  \right), \eqno(5.11)
$$
wherein the (quasi \ - ) optimal values $ N(n) = N_0(n) $ are following: $  N_0(n) = $

$$
 \Ent \left[ \frac{M_{\beta}}{ E_{ \beta, \ 3/2, \ 1/2} \left( \Lambda_{\beta} \ ||K|| \  \Gamma^{3/2}(\beta)  \right)}
\cdot \frac{W_n(\beta) \ ||K||^n \ \Gamma^{1/2}(1 + \beta n) }{\sqrt{n} \ \Gamma^{n/2}(\beta)} \right] + 1. \eqno(5.12)
$$
 Evidently, by the practical using in the case when the value $ N_0(n) $ is sufficiently small, for instance, if $ N_0(n) \le 10, $
then we must  take $ N_0(n):=0; $  see (2.11a). \par

 This imply that the rate of convergence of offered method is equal to $  1/\sqrt{Z}, $ as in the classical Monte Carlo method.\par

 It remains as ordinary to use the classical Central Limit Theorem (CLT) to construct the confidence interval for the solution $  x_{\lambda, \beta}(t). $\par

\vspace{3mm}

\section{ Concluding remarks.}

\vspace{3mm}

{\bf A.  Confidence region for solution in the uniform norm. } \par

\vspace{3mm}

 All the offered (Monte Carlo approximated) solutions, see for example (4.9), $  \hat{x} = \hat{x}(u), \ u \in V $ have a form

$$
\hat{x}(u) =
\hat{x}_Z =\frac{1}{Z} \sum_{i=1}^Z \frac{1}{N(\Delta(i))} \sum_{j=1}^{ N(\Delta(i)) } L_{\Delta(i)} \left(u, \theta^{(j)}_{\Delta(i)} \right), \eqno(6.1)
$$
and are unbiased. \par
 We restrict ourselves for definiteness in this section only the case of the Fredholm's equation; another cases may be considered analogously. \par

 Denote for simplicity

$$
\xi_i(u) = \frac{1}{N(\Delta(i))} \sum_{j=1}^{ N(\Delta(i)) } L_{\Delta(i)} \left(u, \theta^{(j)}_{\Delta(i)} \right) - x(u), \eqno(6.2)
$$
then the random fields $ \{ \xi_i(u) \} $ are continuous, mean zero, identical distributed and

$$
\sqrt{Z} (\hat{x}_Z(u) - x(u)) = \frac{1}{\sqrt{Z}} \sum_{i=1}^Z \xi_i(u). \eqno(6.3)
$$

 In order to build the confidence region in the uniform norm $ ||x|| = \max_{u \in V} |x(u)|  $ for the solution $  x = x(u), $ we need to use the
so - called Central Limit Theorem (CLT) in the space of continuous  functions $  C(V).  $ see
\cite{Ostrovsky1}, \cite{Ostrovsky504}, \cite{Jain1}, \cite{Kozachenko1}, \cite{Ostrovsky302}, \cite{Frolov1}, \cite{Grigorjeva1}.\par
 Namely, if the sequence of random fields $ \xi_i = \xi_i(u), \ u \in V $ satisfies this CLT, then

$$
\lim_{Z \to \infty} {\bf P} \left( \sup_{u \in V} | (\hat{x}_Z(u) - x(u))| > Q   \right) =
{\bf P} \left( \sup_{u \in V} |\zeta(u)| > Q  \right), \eqno(6.4)
$$
where $  \zeta(u)  $  is continuous centered Gaussian random field with at the same  covariation function $  R(u_1,u_2) $ as $ \xi_1(u): $

$$
R(u_1, u_2) = {\bf E} \zeta(u_1) \zeta(u_2) = \Cov(\xi_1(u_1), \xi_1(u_2)).
$$
 Many sufficient conditions for CLT in the Banach space $  C(V) $ may be found in
\cite{Araujo1},   \cite{Dudley1}, \cite{Ledoux1}, \cite{Fortet1},
\cite{Gine1}, \cite{Gine2}, \cite{Heinkel1} etc. CLT in another separable Banach spaces is investigated, e.g. in
\cite{Billingsley2}, \cite{Garling1}, \cite{Gine2}, \cite{Ledoux1}, \cite{Pisier1}, \cite{Ostrovsky401}, \cite{Ostrovsky404},
\cite{Rackauskas1}, \cite{Zinn1}.\par
 The technology of application of the Banach space valued Central Limit Theorem in the parametric Monte Carlo method is described in
 \cite{Frolov1},  \cite{Grigorjeva1}, \cite{Ostrovsky1},  \cite{Ostrovsky303}. \par

\vspace{3mm}

{\bf B. Analogously may be considered the integral equations of a form}

$$
x(t_1,t_2) = f(t_1,t_2) + \int_0^{t_1} ds_1 \int_0^{t_2} ds_2 \cdot K(t_1,t_2, s_1, s_2) \ x(s_1,s_2) \ ds_1 \ ds_2,
$$
with or without weak  singularities.\par

%&


\vspace{4mm}

\begin{thebibliography}{99}

\vspace{4mm}

\bibitem{Ostrovsky1}
{\sc  Ostrovsky E.I.} (1999). {\it Exponential estimations for random Fields and its
Applications, (in Russian).}  Moscow-Obninsk, OINPE.

 \bibitem{Ostrovsky504}
{\sc  Ostrovsky E. and Rogover E.} {\it Exact exponential Bounds for the random Field
Maximum Distribution via the Majorizing Measures (Generic Chaining.) }
 arXiv:0802.0349v1 [math.PR] 4 Feb 2008

\vspace{12mm}

\bibitem{Araujo1}
{\sc Araujo  A., Gine E. }
{\it The central limit theorem for real and Banach valued random variables.}
Wiley, (1980), London, New York.

\bibitem{Billingsley2}
{\sc Billingsley P.} {\it  Convergence of probability measures.}
Wiley, (1968), London, New York.

\bibitem{Dudley1}
{\sc Dudley R.M.} {\it Uniform Central Limit Theorem}. Cambridge University Press, (1999)

\bibitem{Ledoux1}
 {\sc Ledoux M., Talagrand M.} (1991) {\it Probability in Banach Spaces.}
      Springer, Berlin, MR 1102015.

\vspace{12mm}

\bibitem{Fortet1}
{\sc Fortet R. and Mourier E.} {\it Les  fonctions  alratoires comme elements aleatoires dans
les espaces de Banach.} Studia Math., {\bf 15}, (1955), 62-79.

\bibitem{Garling1}
{\sc Garling D.J.H. }
{\it Functional Central Limit Theorems in Banach Spaces.}
 The Annals of Probability, Vol. 4, No. 4 (Aug., 1976), pp. 600-611

\bibitem{Gine1}
{\sc Gine E.} {\it On the Central Limit theorem for sample continuous processes.} Ann.
Probab. (1974), 2, 629-641.

\bibitem{Gine2}
{\sc Gine E., Zinn J.} {\it  Central Limit Theorem and Weak Laws of Large Numbers in certain Banach Spaces.  }
Z. Wahrscheinlichkeitstheory verw. Gebiete. {\bf 62}, (1983), 323-354.

\bibitem{Grothaus1}
{\sc Martin Grothaus,  Florian Jahnert,  Felix Riemann, Jose Luis da Silva.}
{\it Mittag-Leffler Analysis I: Construction and characterization.}
arXiv:1407.8308v1 [math.FA] 31 Jul 2014

\bibitem{Heinkel1}
{\sc Heinkel B.} Measures majorantes et le theoreme de la limite centrale dans
$ C(S). $ Z. Wahrscheinlichkeitstheory. verw. Geb., (1977). 38, 339-351.

\bibitem{Jain1}
{\sc Jain N.C. and Marcus M.B.} {\it Central limit theorem for $C(S)$ valued random
variables.} J. of Funct. Anal., (1975), 19, 216-231.

\bibitem{Kozachenko1}
 {\sc Kozachenko Yu. V., Ostrovsky E.I.} (1985). {\it The Banach Spaces of
      random Variables of subgaussian type.} Theory of Probab. and Math.
      Stat. (in Russian). Kiev, KSU, {\bf 32}, 43-57.

\bibitem{Mittag1}
{\sc Mittag-Leffler G.} {\it Sur la représentation analytique d'une
branche uniforme d'une fonction monogene (cinquieme note).}
Acta Math., {\bf 29(1):},  101 \ - \ 181, 1905.

\bibitem{Ostrovsky301}
{\sc Ostrovsky E., Sirota L.}
{\it CLT for continuous random processes under approximations terms.}
arXiv:1304.0250v1 [math.PR] 31 Mar 2013

\bibitem{Ostrovsky302}
{\sc Ostrovsky E., Sirota L.} {\it  Monte Carlo computation of multiple weak singular integrals of spherical and Volterra's type.  }
arXiv:1405.6344v1 [math.NA] 24 May 2014

\bibitem{Pisier1}
{\sc  Pisier G., Zinn J.} {\it On the limit theorems for random variables with values in the spaces } $ L_p, \ 2  \le p  < \infty. $
 Z. Wahrscheinlichkeitstheorie verw. Gebiete 41, 289 \ - \ 304, (1978).

\bibitem{Rackauskas1}
{\sc Rackauskas A, Suquet Ch.}
{\it Central limit theorems in H\"ölder topologies for Banach space valued random fields.}
Teor. Veroyatnost. i Primenen., 2004, Volume 49, Issue 1, Pages 109 \ – \ 125 (in Russian).

\bibitem{Zinn1}
{\sc Zinn J.  } {\it A Note on the Central Limit Theorem in Banach Spaces.}
 Ann. Probab. Volume 5, Number 2 (1977), 283\ - \ 286.

\vspace{12mm}

\bibitem{Bornemann1}
{\sc Bornemann Folkmar} {\it  On the numerical evaluation of Fredholm determinants.}
 arXiv:0804.2543v1 [math.NA] 16 Apr 2008

\bibitem{Fredholm1}
{\sc Fredholm I.} {\it Sur une classe d'equations  fonctionelles. } (1903), Acta Math., {\bf 27,} 365-390.

\bibitem{Frolov1}
{\sc Frolov A.S., Tchentzov N.N. } {\it On the calculation by the Monte-Carlo
method definite integrals depending on the parameters. } Journal of Computational
Mathematics and Mathematical Physics, (1962), V. 2, Issue 4, p. 714-718 (in
Russian).

\bibitem{Dunford1}
{\sc Dunford N., Schwartz J.T.  } Linear operators. Part 1, General Theory, (1958),  Interscience Publishers,
New York, London.

\bibitem{Grigorjeva1}
{\sc Grigorjeva M.L., Ostrovsky E.I.} {\it Calculation of Integrals on discontinuous
Functions by means of depending trials method.} Journal of Computational
Mathematics and Mathematical Physics, (1996), V. 36, Issue 12, p. 28-39 (in
Russian).

\bibitem{Ostrovsky303}
{\sc Ostrovsky E., Sirota L.} {\it Monte-Carlo method for multiple parametric integrals
calculation and solving of linear integral Fredholm equations of a second
kind, with confidence regions in uniform norm.}
 arXiv:1101.5381v1 [math.FA] 27 Jan 2011

\vspace{12mm}

\bibitem{Ostrovsky401}
{\sc Ostrovsky E., Sirota L.} {\it  Central Limit Theorem and exponential tail estimations in mixed (anisotropic) Lebesgue spaces. }
arXiv:1308.5606v1 [math.PR] 26 Aug 2013

\bibitem{Ostrovsky404}
{\sc Ostrovsky E., Sirota L.} {\it  Central Limit Theorem anf exponential tail estimates in hybrid Lebesgue-continuous spaces.  }
arXiv:1309.2344v1 [math.PR] 9 Sep 2013

\vspace{4mm}

\end{thebibliography}
\end{document}